\documentclass[12pt,reqno]{amsart}
\usepackage{amsmath}
\usepackage{amssymb,amscd}
\usepackage{multicol}
\usepackage{hyperref}
\usepackage[all,cmtip]{xy}
\numberwithin{equation}{section}

\newtheorem{theorem}{Theorem}[section]
\newtheorem{corollary}[theorem]{Corollary}

\theoremstyle{definition}
\newtheorem{defn}[theorem]{Definition}

\def \({\left(}
\def \){\right)}
\def \<{\langle}
\def \>{\rangle}

\newcommand{\End}{\mbox{End}}

\begin{document}

\title[Toric co-Higgs bundles on toric varieties]{Toric co-Higgs bundles on toric 
varieties}

\author[I. Biswas]{Indranil Biswas}

\address{School of Mathematics, Tata Institute of Fundamental Research,
Homi Bhabha Road, Mumbai 400005, India}

\email{indranil@math.tifr.res.in}

\author[A. Dey]{Arijit Dey}

\address{Department of Mathematics, Indian Institute of Technology-Madras, Chennai 600036, India}

\email{arijit@iitm.ac.in}

\author[M. Poddar]{Mainak Poddar}

\address{Mathematics Department, Indian Institute of Science Education and Research, Pune 411008, India}

\email{mainak@iiserpune.ac.in}

\author[S. Rayan]{Steven Rayan}

\address{Centre for Quantum Topology and Its Applications (quanTA) and Department of Mathematics \& 
Statistics, University of Saskatchewan, McLean Hall, 106 Wiggins Road, Saskatoon, SK S7N 5E6, 
Canada}

\email{rayan@math.usask.ca}



\date{July 27, 2020}

\begin{abstract} 
Starting from the data of a nonsingular complex projective toric variety, we define an associated notion of \emph{toric co-Higgs 
bundle}. We provide a Lie-theoretic classification of these objects by studying the interaction between Klyachko's fan filtration and the fiber of the co-Higgs bundle at 
a closed point in the open orbit of the torus action. This can be interpreted, under certain conditions, as the construction of a 
coarse moduli scheme of toric co-Higgs bundles of any rank and with any total equivariant Chern class.
\end{abstract}

\maketitle

\section{Introduction\\}

\vspace{10pt}

We begin with an algebraic or, equivalently, holomorphic vector bundle $V$ over a nonsingular complex projective 
variety $X$ with tangent bundle $TX$. Then, a \emph{co-Higgs field} for $V$ is a 
holomorphic section $\phi$ of the twisted endomorphism bundle $\End(V) \otimes TX$, subject to the integrability condition that the quadratic section $\phi\otimes\phi$ is symmetric --- that is, that the 
section $\phi \wedge \phi$ of $\End(V) \otimes \Lambda^2 TX$ vanishes identically. A pair 
$(V,\,\phi)$ satisfying the above conditions is referred to as a \emph{co-Higgs bundle}. Co-Higgs 
bundles were introduced simultaneously by Hitchin \cite{Hit2} and the fourth-named author 
\cite{Ray1} in the context of generalized complex geometry. The name \emph{co-Higgs} speaks to a 
duality with Higgs bundles in the sense of Hitchin \cite{Hit0,Hit1} and Simpson \cite{Sim1}, where 
the Higgs fields are $T^*X$-valued.\\

Co-Higgs bundles have been classified and/or constructed on $\mathbb P^1$ \cite{Ray2,BGHR}, 
$\mathbb P^2$ \cite{Ray3}, $\mathbb P^1\times\mathbb P^1$ \cite{VC}, and logarithmic curves 
\cite{BH}, for example. Over singular varieties, they have been used to some effect towards establishing inequalities related to vector-valued modular forms \cite{FR}. At the same time, there are ``no-go'' theorems for the existence of nontrivial co-Higgs bundles in some instances, such as over the moduli space of stable bundles on a nonsingular complex curve of 
genus at least $2$ \cite{BR1}. Recently, the first-named and fourth-named authors classified 
homogeneous co-Higgs bundles on Hermitian symmetric spaces \cite{BR2}. The goal of the present work 
is to extend this to toric varieties. Accordingly, we define below a natural notion of \emph{toric co-Higgs bundle}.\\

We classify toric co-Higgs structures for a fixed toric bundle $V$ using Klyachko's seminal work on classification of toric vector bundles \cite{Kly}. We recall from \cite[Theorem 2.2.1]{Kly} that the category of toric bundles on $X$ with 
equivariant homomorphisms is equivalent to the category of compatible $\Sigma$-filtered vector 
spaces, where $\Sigma$ is the fan of $X$. Our strategy for classifying toric co-Higgs bundles 
$(V,\,\phi)$ is thus to reduce the data of $(V,\,\phi)$ to that of a tuple of commuting 
$\Sigma$-filtered endomorphisms of the fiber $V_{x_0}$, where $x_0$ is a closed point in the open 
orbit of $T$ in $X$. This is the content of Theorem \ref{mainthm}, which is the main theorem in this note.\\

Similar to the symmetric space case \cite{BR1}, the resulting classification is Lie-theoretic in nature 
and admits an interpretation as a moduli construction. Subject to certain conditions (namely, the freeness of a certain group action), we identify a scheme 
that parametrizes toric co-Higgs bundles of fixed rank and total equivariant Euler characteristic in the sense of \cite{Pay1}. This
scheme fibers over an associated moduli scheme of toric bundles constructed in \cite{Pay2}.\\

We also wish to point out similar work in \cite{AW}, completed independently and at roughly the same time.\\

\section{Set-up and examples\\}

\vspace{10pt}

\subsection{Basic notions} Throughout, $X$ is a nonsingular complex projective variety. Assume that $X$ admits an algebraic (equivalently, holomorphic) action of a complex torus $T\,\cong\, (\mathbb C^*)^n$ so that it is a \emph{toric variety} in the sense of \cite{Ful}. Furthermore, fix a holomorphic vector bundle $V$ and suppose that it admits a 
lift of the action of $T$ from $X$ which is fiber-wise linear. In other words, $V$ is equipped with 
the structure of a $T$-equivariant vector bundle.
We will refer to $V$ simply as a \emph{toric bundle}. \\

There is subsequently an induced action of $T$ on the vector space of global holomorphic sections of $V$:
\begin{equation}\label{eqaction}
(t \cdot s)(x) \,=\, t s( t^{-1} x)
\end{equation}
for all $s\,\in\, H^0(X,\, V)$ and $t\,\in\, T$. A section $s\,\in\, H^0(X,\, V)$ is 
said to be \emph{semi-invariant} if there exists a character $\chi$ of $T$ such that 
\begin{equation}\label{eqchi}
t\cdot s \,=\, \chi(t) s
\end{equation}
for all $t\,\in\, T$. A semi-invariant section $s$ is said to be \emph{invariant}
if the associated character $\chi(t)$ is trivial, meaning 
 \begin{equation}\label{eqinvsec}
 t\cdot s \,=\, s
 \end{equation}
for all $t\,\in\, T$. Combining \eqref{eqaction} and \eqref{eqinvsec}, we conclude that 
\begin{equation}\label{eqchi2}
 t s( t^{-1} x) \,=\, s(x)
\end{equation}
for any invariant section $s$.\\

\subsection{$T$-equivariant structures and toric co-Higgs bundles}

Now, a $T$-equivariant structure on $V$ induces a $T$-equivariant structure on $\End(V)$ in the following way. Given an
element $\psi$ in the fiber $(\End(V))_x\,=\,\End(V_x)$ over $x\,\in\, X$, we
define $t\psi$ in $\End(V_{tx})$ by
\begin{equation}\label{eqbullet}
(t \psi) (v) \,=\, t(\psi(t^{-1}v))
\end{equation}
for every $v \,\in\, V_{tx}$ and $t\,\in\, T$. Then, by \eqref{eqchi2}, a holomorphic
section $\phi$ of $\End(V)$ is invariant if and only if 
\begin{equation}\label{eqequiv}
t\phi(t^{-1} x) \,=\, \phi(x)
\end{equation}
for all $x \,\in\, X$ and $t\,\in \,T$. By \eqref{eqbullet}, equation \eqref{eqequiv} is equivalent to 
$$ t \circ \phi(t^{-1}x) = \phi(x) \circ t $$ for all $x \,\in\, X$ and $t\,\in \,T$.
In other words, $\phi$ is an invariant section of $\End(V)$
if and only if $\phi$ is a $T$-equivariant endomorphism of $V$.\\

At the same time, the tangent bundle $TX$ has a natural $T$-equivariant structure that is simply 
the linearization of the $T$-action on $X$. Together with the above $T$-action on $\End(V)$, we 
have an induced $T$-equivariant structure on the twisted bundle $\End(V)\otimes TX$. This allows 
us to formulate the following:\\
 
\begin{defn}
A \emph{toric co-Higgs bundle} on a toric variety $X$ is a pair $(V,\,\phi)$, where $V\,\longrightarrow\, X$ is
a toric bundle and $\phi\, \in\, H^0(X,\, \End(V) 
\otimes TX)$ is an invariant co-Higgs field.\\
\end{defn}

\subsection{Examples}\label{example1}

The tangent bundle $TX$ of a nonsingular projective toric variety $X$ 
always admits a nonzero invariant holomorphic section $s$. Let us take the tensor product of such 
a holomorphic vector field with the identity homomorphism $1$ of \emph{any} toric bundle $V$. Call this product $\phi$. It 
follows immediately that
$$\phi\wedge\phi\,=\,(s\otimes 1)\wedge(s\otimes1)\,=\,[s,\,s]\otimes(1\wedge1)\,=\,0\, .$$
In other words, $\phi$ is a nontrivial invariant co-Higgs field on $V$, and hence every toric vector bundle on any toric variety is equipped with a family of invariant co-Higgs fields induced by the invariant holomorphic vector fields on $X$.\\

Another example, similar in spirit but for a specific bundle, is given by choosing $V\,=\,TX\oplus\mathcal O_X$ with its natural toric structure, where 
$\mathcal O_X$ is the structure sheaf of $X$. We can equip this $V$ with the co-Higgs field
$$\phi\,=\,\left(\begin{array}{cc}0 & 1\\0 & 0\end{array}\right)\, ,$$
where $1$ is interpreted as the identity 
morphism on $TX$, which serves as the part of the co-Higgs structure that acts as $TX\,\longrightarrow\, \mathcal 
O_X\otimes TX$. Because $\phi$ is built from just the identity map, the co-Higgs field automatically has the required invariance. It satisfies the vanishing condition $\phi\wedge\phi\,=\,0$ since, as a matrix, $\phi$ is 
nilpotent. More generally, this example is present on any complex variety --- it is the so-called \emph{canonical co-Higgs bundle} --- and it is 
discussed in some detail from the point of view of slope stability and deformation theory in Chapter 6 of \cite{Ray1}.\\ 

For contrast, we conclude with an example of a co-Higgs structure on a vector bundle that is not directly related to toric
geometry. We briefly recall the definition of the jet bundle of a holomorphic 
vector bundle $E$ on a projective variety $Y$. Let 
\[ p_j:\ Y \times Y \longrightarrow Y, \,\, j\,=\,1,2 \]
be the natural projection to the $j$-th factor. Let $\Delta$ be the diagonal 
of $Y\times Y$. For any non-negative integer $k$, consider the $k$-th jet bundle of $E$ defined as 
\[ J^k(E) \,=\, p_{1 \ast}(p_2^{\ast}(E)/(p_2^{\ast}E \otimes \mathcal O_{Y\times Y}(-(k+1)\Delta) \,\longrightarrow\, Y
\]
One can check that the higher direct image of the above right hand side vanishes, which makes $J^k(E)$ a vector bundle. 
 We have the following exact sequence 
\[
\xymatrix{
	0 \ar[r] & E \otimes \Omega_{Y} \ar[r]^{f} & J^1(E) \ar[r]^{g} & J^{0}(E) 
	\ar[r]& 0
} 
\]
 The composition homomorphism $(g \otimes 1_{\Omega_{Y}}) \circ f$ 
defines a nonzero co-Higgs structure on $J^1(E)$. Even when $Y$ itself is toric, this gives examples
of non-toric co-Higgs bundles by taking $E$ to be non-toric.  For further examples of non-toric co-Higgs bundles on toric varieties, we refer the reader to \cite{AW}.\\

In the next section, we classify toric co-Higgs structures on a fixed toric bundle.\\

\section{Toric co-Higgs bundles and the Klyachko fan filtration\\}

\vspace{10pt}
 
\subsection{$\Sigma$-filtrations}

Let $X$ be a nonsingular complex projective toric variety equipped with an action of $T$ (cf. \cite{Ful}). Let $\Sigma$ denote the fan of $X$, and let $\Sigma(d)$ be the set of 
$d$-dimensional cones in $\Sigma$. For any toric bundle 
$V$ on $X$, Klyachko \cite{Kly} constructed a compatible family of full filtrations of 
decreasing subspaces of the fiber $E \,=\, V_{x_0}$, where $x_0$ is a closed point in the open 
orbit of $T$ in $X$. We will subsequently refer to $x_0$ simply as a ``closed point for $T$''. The 
family is indexed by the $T$-invariant divisors or equivalently by $|\Sigma(1)|$, the set of integral generators of the one-dimensional cones of $\Sigma$. In other words, we have a family of filtrations
$$\{E^{\rho}(i) \,\mid\, \rho \,\in \,|\Sigma(1)|,\ i \,\in\, \mathbb{Z}\}\, .$$
 Note that ``decreasing'' means $$E^{\rho}(i+1) 
\,\subseteq\, E^{\rho}(i)$$ for all $i$. For brevity, such a family 
of filtrations will be called a \emph{$\Sigma$-filtration}. The compatibility condition mentioned above 
refers to the existence of cone-wise $T$-module structures on $E$ giving rise to the 
$\Sigma$-filtration. We give below a short review of the construction of the $\Sigma$-filtration
in the case of a projective toric variety. In this case, the description is slightly simpler as it suffices to consider only the cones of dimension $n$ in the fan.\\

Let $M$ and $N$ denote the dual lattices of characters and co-characters of $T$, respectively. 
Let $\sigma$ be a maximal cone of $\Sigma$ and $X_{\sigma}$ the corresponding affine toric variety.
Denote by $V_{\sigma}$ the restriction of $V $ to $X_{\sigma}$.
It was shown by Klyachko \cite[Proposition 2.1.1]{Kly} that there exists
 a framing (which is not unique) of $V_{\sigma}$ by semi-invariant sections. 
 Fix such a framing $(s_1, \ldots, s_r)$.
 Let $S_{\sigma}$ be the $T$-submodule of
 $H^{0}(X_{\sigma},\, V_{\sigma})$ generated by the semi-invariant sections $s_1, \ldots, s_r$.
Evaluation at $x_0$
gives an isomorphism of vector spaces $ev_{x_0}\,:\, S_{\sigma}\,\longrightarrow \, E$.
 This isomorphism induces a $T$--module structure on $E$, or equivalently, a decomposition
\begin{equation}\label{dec1}
E \,=\, \bigoplus_{u \in M} E^{\sigma}_u\, .
\end{equation} 
The decompositions \eqref{dec1} may depend on the choice of the semi-invariant framing of $V_{\sigma}$.
However, Klyachko showed that for each $\rho \in |\Sigma(1)|$, the subspaces
\begin{equation}\label{filt}
E^{\rho}(i) \,:=\, \bigoplus_{u\in M, u(\rho) \ge i} E^{\sigma}_{u}\, ,\quad {\rm where} \; \sigma \in \Sigma(n) \; {\rm is \; such \; that}\;
\rho \,\in \,|\Sigma(1)| \bigcap \sigma \, , \end{equation}
 are independent of the choice of
$\sigma$ containing $\rho$ and the choice of the framing. \\

 A morphism of compatible 
$\Sigma$-filtered vector spaces $\{E^{\rho}(i)\}$ and $\{F^{\rho}(i)\}$ is a vector space map 
$\phi\,: E \,\longrightarrow\, F$ such that $\phi(E^{\rho}(i))\,\subseteq\, F^{\rho}(i)$ for all 
$\rho$ and $i$. We call such a morphism a \emph{filtered linear map} of $\Sigma$-filtered vector 
spaces. Proposition 2.1.1(iii) of \cite{Kly} proves that equivariant morphisms between two toric vector bundles over $X$ correspond to filtered linear maps between their $\Sigma$-filtrations. Note that any linear map between the fibers at $x_0$ extend uniquely to a $T$-equivariant morphism of vector bundles
over the open $T$-orbit. The regularity of such an extension over the boundary of the open orbit is naturally and intimately related to the weights of the $T$-action and the defining inequalities of the $\Sigma$-filtrations.\\

The discussion above can be summarized as following celebrated theorem of Klyachko. \\

\begin{theorem}[{\cite[Theorem 2.2.1]{Kly}}]\label{Kly}
The category of toric vector bundle over a toric variety $X$ with fan $\Sigma$ is equivalent to the category of complex vector spaces $E$ with a family of decreasing filtrations $$\{E^{\rho}(i) \,\mid\, \rho \,\in \,|\Sigma(1)|,\ i \,\in\, \mathbb{Z}\}\, ,$$ which satisfy following compatibility condition : 

For each $\sigma \in \Sigma(n)$, there exists a $M$-grading $E\,=\, \bigoplus_{u \in M} E^{\sigma}_u \, $ for which 
$$
E^{\rho}(i) \,:=\, \bigoplus_{u\in M, u(\rho) \ge i} E^{\sigma}_{u}\, ,\quad {\rm where} \; \sigma \; {\rm is \; such \; that}\;
\rho \,\in \,|\Sigma(1)| \bigcap \sigma. \, 
$$
\end{theorem}
\vspace{.5cm}

\subsection{The main theorem}

Having built up necessary language around $\Sigma$-filtrations, we now pose the main theorem of this article.

\begin{theorem}\label{mainthm}
Let $X$ be a nonsingular complex projective toric variety equipped with an action of $T\,\cong\,(\mathbb C^*)^n$, let $V$ by any toric bundle on $X$, and let $x_0\, \in\, X$ be a closed point for $T$. Then, there is a $1:1$ correspondence between invariant co-Higgs fields $\phi$ and $n$-tuples of pairwise-commuting filtered linear maps of $E\,= \,V_{x_0}$ that respect the Klyachko $\Sigma$-filtration.\\
\end{theorem}

Before proving Theorem \ref{mainthm}, we need to understand the integrability condition $\phi\wedge\phi\,=\,0$ locally. In \cite{Ray1,Hit2,Ray3}, a local criterion for the 
vanishing of $\phi \wedge \phi$ is identified. Suppose that $\{z_1,\, \dots,\, z_n\}$ is a 
holomorphic coordinate system on an affine chart $U$ in a nonsingular variety $X$. We can write
$$
\phi|_U \,=\, \sum_{i =1}^n \phi_i \frac{\partial}{\partial z_i},
$$where each $\phi_i \in H^0(U,\,\End(V))$. Then 
\begin{equation}\label{eqcommutator}\phi \wedge \phi = 0 \; {\rm on}\; U \iff [\phi_i, \phi_j ] = 0 \; {\rm on}\; U \; \forall \, 1 \le i,j \le n\,.
\end{equation}

\vspace{10pt}

With this observation in hand, we can proceed with the proof of the main theorem.\\

\begin{proof}[{Proof of Theorem \ref{mainthm}}]
Let $(t_1, \,\dots,\, t_n)$ be coordinates on $T$ corresponding to an integral basis of 
$\mbox{Lie}(T)$. We identify these with coordinates on the open dense $T$-orbit in $X$; this 
should not cause any confusion. To conform to \eqref{eqaction}, we consider the action of $T$ on 
$\mathcal{O}_T$ (and $\mathcal{O}_X$) given by 
\begin{equation}\label{function} (t\cdot f)(x) := f(t^{-1}x)\,.
\end{equation}
Observe that the character $t_i$, and the derivation $\frac{\partial}{\partial t_i}$, have weights $t_i^{-1}$ and $t_i$ respectively under this action.
The corresponding invariant vector fields $t_i \frac{\partial}{\partial t_i}$ are naturally $T$-invariant on the open orbit in $X$. By \cite[Theorem 3.1]{BDP}, these vector fields admit 
$T$-invariant holomorphic extensions to the whole of $X$. In fact, by Lie theory, any $T$-invariant vector field on $X$ is a complex linear combination of these fields.\\

Now, let $A_1,\, \dots,\, A_n$ be pairwise-commuting linear endomorphisms of $V_{x_0}$ that 
respect the $\Sigma$-filtration. Then by Klyachko's theorems \cite[Proposition 2.1.1, Theorem 2.2.1]{Kly}, these define $T$-equivariant 
endomorphisms $\phi_1, \dots, \phi_n $ of $V$ such that $\phi_j(x_0)\,=\, A_j$. Therefore, each 
$\phi_j$ is an invariant section of $\End(V)$. Applying equation \eqref{eqchi2} to $\phi_j$ we 
see that
$$
\{t \phi_j(t^{-1} x_0) \,= \, \phi_j(x_0)\} \implies \{\phi_j (t^{-1} x_0) \,= \,t^{-1} \phi_j(x_0)\}$$ for 
all $t\,\in\, T$. It then follows from \eqref{eqbullet} that the $\phi_j(t^{-1} x_0)$'s commute with 
each other for every $t \,\in\, T$. In other words, they commute mutually on the open dense $T$-orbit 
in $X$. Therefore, by continuity, the $\phi_j$'s commute on entire $X$.\\

Next, we define $ \phi \,\in\, H^0 (X,\,\End(V) \otimes TX)$ by 
\begin{equation}\label{eqphi}
\phi \,=\, \sum_{j} \phi_j \, t_j\, \frac{\partial}{\partial t_j} \,.
\end{equation} 
Consider any affine toric chart on $X$ with coordinates $(z_1,\, \dots,\, z_n)$. By \cite[Lemma 3.1]{BDP} we have 
\begin{equation}\label{eqcoc} 
t_j\, \frac{\partial}{\partial t_j} \,=\, \sum_k c_{jk}(z_1, \dots, z_n) \frac{\partial}{\partial z_k}, \, \end{equation} 
where the $c_{jk}$'s are holomorphic functions.
Substituting \eqref{eqcoc} in \eqref{eqphi}, we have the following representation of 
$\phi$ in the $(z_1,\, \dots,\, z_n)$ coordinates:
$$ \phi \,=\, \sum_k \psi_k \frac{\partial}{\partial z_k}\, ,$$
where $\psi_k \,=\, \sum_j c_{jk}(z_1, \dots, z_n ) \phi_j$.
Since the $\phi_j$'s commute and the $c_{jk}$'s are scalars, the $\psi_k$'s also mutually commute. 
Thus by \eqref{eqcommutator}, $\phi$ defines a co-Higgs structure on $V$. Hence, given a tuple
$(A_1,\, \dots,\, A_n)$ of commuting filtered endomorphisms of $V$, we obtain an
equivariant co-Higgs structure on $V$.\\

In the other direction, given any equivariant co-Higgs structure $\phi$ on $V$, we may write 
$\phi$, on the open orbit, as in \eqref{eqphi}. We use the fact any torus-equivariant vector bundle 
is trivial over the open orbit. As the vector fields $ t_j \frac{\partial}{\partial t_j} $ are 
$T$-invariant, the $\phi_j$'s are also $T$-invariant. Moreover, as the open orbit is an 
affine toric chart, and $\phi$ is a co-Higgs field, the $\phi_j$'s commute mutually by 
\eqref{eqcommutator}. Then we define $A_j\,=\,\phi_j(x_0)$. As $\phi_j$ is a $T$-equivariant 
 endomorphism of $V$, the endomorphisms $A_j$ respect the Klyachko $\Sigma$-filtration.\\

It is now straightforward to check that the above association is a bijection.
\end{proof}

\subsection{Examples for the theorem}

As an example, we take $V$ to be the tangent bundle $T X$. We will describe the toric co-Higgs fields on $V$ for certain toric varieties $X$, and also some co-Higgs fields that are not invariant but only semi-invariant. Let 
$\rho_1, \ldots, \rho_m$ denote the elements of $|\Sigma(1)|$. We may assume without loss of generality that $\rho_1, \ldots, \rho_n$ form an integral basis of the co-character lattice of $T$.
We choose the closed point $x_0$ in the open orbit to correspond to the identity element of $T$. Then
 the underlying vector space $E=V_{x_0}$ of the Klyachko $\Sigma$-filtration may be identified with 
$$\mathbb{C}\left\langle \frac{\partial}{\partial t_1},\, \ldots ,\, \frac{\partial}{\partial t_n} \right\rangle \cong 
\mathbb{C}\langle \rho_1,\,\ldots,\,\rho_n \rangle\,.$$ Recall that under the action \eqref{function} the derivation
$\frac{\partial}{\partial t_j}$ has weight $t_j$. It then follows from \eqref{filt} that the Klyachko filtration on $T\mathbb{P}^n$ (cf. \cite[p. 350]{Kly}) is given by 
$$ E^{\rho_j}(i)\,=\, \left\{ \begin{array}{ll}
0 & {\rm if \,} i > 1 \\
\mathbb{C}\langle \rho_j \rangle & {\rm if \,} i = 1\\
E & {\rm if \,} i \le 0 \,.\end{array}\right. $$ 
Now, let $X\,=\, \mathbb{P}^n$.  Then, we have $m \,=\, n+1$ and $\rho_{n+1}\,=\, -\sum_{j=1}^n \rho_j$. In this case, the only filtered endomorphisms of $E$ are constant multiples of the identity map. Thus, any toric co-Higgs field on $T\mathbb{P}^n$ is of the form $s \otimes 1$ where $s= \sum a_j t_j \frac{\partial }{\partial t_j}$ for some $a_j \in \mathbb{C}$.  However, $T\mathbb{P}^n$ admists others co-Higgs fields. For instance, each $\frac{\partial}{\partial t_j}$ is a global semi-invariant section of $T\mathbb{P}^n$ with weight $t_j$. Thus, $\frac{\partial}{\partial t_j} \otimes 1$ is a semi-invariant co-Higgs field on the bundle $T\mathbb{P}^n$, which is not toric according to our definition. Moreover, 
$\left(\frac{\partial}{\partial t_1} + \frac{\partial}{\partial t_2} \right) \otimes 1$ is a co-Higgs field 
on $T\mathbb{P}^n$, for $n\,\ge\, 2$, which is not even semi-invariant.\\

More generally, we obtain a larger class of filtered endomorphisms of $E$ when $X$ is the product of projective spaces. For instance, let $X = \mathbb{P}^1 \times \mathbb{P}^1$. In this case, 
$|\Sigma(1)| = \{\rho_1, \ldots, \rho_4 \}$ where $\rho_1\,=\, (1,\,0)\, =\, -\rho_3$ and 
$\rho_2 \,=\,(0,\,1)\, =\, -\rho_4 $. Thus, $E \cong \mathbb{C} \langle \rho_1\rangle \oplus \mathbb{C} \langle \rho_2 \rangle$ and filtered endomorphisms of $E$ are given by diagonal matrices with respect to this decomposition. Since any two diagonal matrices commute, a toric co-Higgs field on $TX$ is of the form $\phi_1 t_1 \frac{\partial}{\partial t_1} +\phi_2 t_2 \frac{\partial}{\partial t_2}$ where 
$\phi_1$ and $\phi_2$ are arbitrary filtered endomorphisms of $E$. 

\section{Existence of a moduli scheme\\}

\vspace{10pt}

We fix a toric bundle $V$ on a non-singular toric variety $X$ of dimension $n$ with fan $\Sigma$. Let $x_0$ be a closed point of $T$ and put $E=V_{x_0}$ as in the preceding section.
We use $H^{\rho}$ to refer to the parabolic subgroup of $\mbox{GL}(E)$ that preserves the filtration $E^{\rho}$ 
on $E$. Then, the group of endomorphisms of the $\Sigma$-filtration $\{ E^{\rho} (i)\}$ coincides 
with the group $H_V \,:=\, \bigcap_{\rho} H^{\rho}$. Notice that the group $H_V$ contains the center of 
$\mbox{GL}(E)$. Denote by $H_V[n]$ the set of $n$-tuples of pairwise-commuting
elements of the group $H_V$. Now, Theorem \ref{mainthm} can be recast as:\\

\begin{corollary}
If $X$ is a nonsingular complex projective toric variety and $V$ is any toric bundle on $X$, then invariant co-Higgs fields $\phi$ for $V$ are in $1:1$ correspondence with elements of $H_V[n]$.\\
\end{corollary}

Now, we wish to consider all isomorphism classes $[V]$ of toric bundles on $X$ having fixed rank $r$ and fixed $T$-equivariant Chern classes. These classes can be defined explicitly within the equivariant Chow cohomology ring of $X$ in terms of Klyachko's filtration as per \cite{Pay1}. As per \cite{Pay2}, let $\mathcal V^{fr}_X(r,\psi)$ be the fine moduli space of rank-$r$ toric vector bundles framed at $x_0$ and with total equivariant Chern class $\psi$. It is then a result of Payne \cite[Corollary 3.11]{Pay2} that if $\mbox{PGL}(r)$ acts freely on $\mathcal V^{fr}_X(r,\,\psi)$, then there exists an associated coarse moduli scheme $\mathcal V_X(r,\psi)$ of toric bundles on $X$ with that total equivariant Chern class. In light of this, our result gives rise to:\\

\begin{corollary} When the group $\mbox{\emph{PGL}}(r)$ acts freely on $\mathcal V^{fr}_X(r,\,\psi)$, there exists a scheme$${\mathcal C}_X(r,\psi)\,\stackrel{\pi}{\longrightarrow}\,{\mathcal V}_X(r,\psi),$$with fibers $\pi^{-1}([V])\, \cong \,H_V[n]$, that can be identified with a quasiprojective coarse moduli scheme of toric co-Higgs 
bundles on $X$ of fixed rank $r$ and total equivariant Chern class $\psi$.\\\end{corollary}

We take the opportunity now to mention a few natural questions for further study. First, assuming the fibration $\pi$ exists, when is it flat? 
Recall that flat morphisms, in reasonable circumstances, have strong topological properties. 
For instance, the fibers of a surjective, faithfully flat morphism of irreducible varieties will have the expected dimension, which is the difference in the dimensions of the ambient schemes.\\

Second, does ${\mathcal C}_X(r,\,\psi)$ inherit arbitrarily bad singularities from ${\mathcal V}_X(r,\psi)$, as per the ``Murphy's Law'' 
for toric bundles in \cite[Section 4]{Pay2}, and under which condition(s) does it becomes smooth? It is known that in the rank-$2$ case that the framed moduli
$\mathcal V^{fr}_X(2,\psi)$ is smooth.  It is natural to ask whether this holds for co-Higgs case.\\

Moreover, it would be desirable to understand the relationship of this construction to either Mumford-Takemoto or Gieseker stability for co-Higgs bundles in general. In particular, Simpson's moduli space of $\Lambda$-modules \cite{Sim2}, where $\Lambda$ is a coherent sheaf of $\mathcal O_X$-modules, produces a moduli space of Gieseker-stable 
coherent co-Higgs sheaves on $X$ when $\Lambda\,=\,\mbox{Sym}^\bullet(T^*X)$ (cf. \cite[Section 2]{Ray3} for further details on this 
correspondence in the co-Higgs setting). Interpreting the variation of toric structures on $V$ with regards to the moduli problem for 
$\Lambda$-modules is an interesting direction for further exploration. Finally, classifying toric co-Higgs bundles over singular toric 
varieties is a natural problem to explore.

\section*{Acknowledgements}

\vspace{10pt}

We thank Kael Dixon for useful conversations and the anonymous referee for helpful suggestions. The first-named author is partially supported by a J. C. Bose Fellowship, and 
the School of Mathematics, TIFR, is supported by 12-R$\&$D-TFR-5.01-0500. The second-named author is supported in part by 
a SERB MATRICS Grant, MTR/2017/001008/MS.
The research of the third-named author is supported in part by a SERB MATRICS Grant, 
MTR/2019/001613. The fourth-named author is partially supported by a Natural Sciences and Engineering Research Council of Canada (NSERC) Discovery Grant.

\end{document}